\documentclass[a4paper]{article}
\usepackage[left=4cm, right=4cm, bottom=3.5cm, top=3.5cm]{geometry}
\usepackage{graphicx}

\providecommand{\keywords}[1]{\textit{Keywords: } #1}

\usepackage{hyperref}
\usepackage{amsfonts, amsbsy, amssymb, amsmath}

\begin{document}

\markboth{Kraj{\v{n}}{\'a}k and Wiggins}{Dynamics of the Morse Oscillator: Analytical Expressions for Trajectories, Action-Angle Variables, and Chaotic Dynamics}

\title{Dynamics of the Morse Oscillator: Analytical Expressions for Trajectories, Action-Angle Variables, and Chaotic Dynamics}

\author{Vladim{\'i}r Kraj{\v{n}}{\'a}k\footnote{School of Mathematics,
University of Bristol, Bristol BS8 1TW, United Kingdom}, Stephen Wiggins\footnotemark[1]}

\maketitle

\begin{abstract}
We consider the one degree-of-freedom Hamiltonian system defined by the Morse potential energy function (the ''Morse oscillator''). We use the geometry of the level sets to construct explicit expressions for the trajectories as a function of time, their period for the bounded trajectories, and action-angle variables. We use these trajectories to prove sufficient conditions for chaotic dynamics, in the sense of Smale horseshoes, for the time-periodically perturbed Morse oscillator using a Melnikov type approach.

\keywords{Morse oscillator, action-angle variables, homoclinic orbit, Melnikov function, chaos.}
\end{abstract}

\section{Introduction}
\label{sec:intro}

The potential energy function derived by P. M. Morse is truly a ''workhorse''  potential energy function in theoretical chemistry \cite{morse1929diatomic}.  Originally it was devised to describe the intermolecular force between the two atoms in a diatomic molecule. It has the functional form:

\begin{equation}
V(q) =  D \left( 1-e^{-\alpha q} \right)^2,
\label{eq:mpot}
\end{equation}

\noindent
where $q$ represents the distance between the two atoms, $D >0$ represents the depth of the potential well (defined relative to the dissociated atoms), and $\alpha >0$ controls the width of the potential well ($\alpha$ small corresponds to a ''wide'' well, $\alpha$ large corresponds to a narrow well). 

The Morse potential defines a one degree-of-freedom Hamiltonian system, i.e. the phase space is two dimensional described by coordinates $(q, p)$, where $p$ is the momentum conjugate to the position variable $q$.  The Hamiltonian has the form of the sum of the  kinetic energy and the potential energy (the Morse potential).  The Hamiltonian system is integrable, and all trajectories lie on the level sets of the Hamiltonian function.  The level sets can be used to derive integrals that give (time) parametrizations of trajectories. Regions of closed (bounded) trajectories can be used to construct special coordinates--{\em action-angle coordinates}, where the angle denotes a particular location on the closed level set and the action is the area enclosed by a closed level set (divided by $2 \pi$). The transformation to action-angle coordinates for integrable Hamiltonan systems is a standard topic in good classical mechanics textbooks, see, e.g. \cite{landau1960mechanics, arnold2013mathematical}. The transformation preserves the Hamiltonian nature of the system, i.e. it is a canonical transformation, and therefore the standard approach to constructing such transformations is through the use of generating functions. However, for one degree-of-freedom time-independent Hamiltonian systems (such as the one described by the Morse potential)  there is a simpler approach to generating the action-angle transformation that uses the geometry of the closed level set of the Hamiltonian function and the explicit (time) parametrization of the trajectories that can  be obtained (in principle, if the necessary integrals can be explicitly computed)  for one degree-of-freedom Hamiltonian systems. The approach is inspired by the seminal paper of Melnikov \cite{melnikov1963vk}. This approach was developed in detail in \cite{wiggins1990introduction} (the 1990 edition, {\em not} the 2003 edition)
 and is also described in \cite{mezic1994integrability}. This is the approach that we will follow here.

 Action-angle variables are important in Hamiltonian mechanics for a number of reasons. From the point of view of classical mechanics they are the coordinate system used for the development of the Kolmogorov-Arnold-Moser and Nekhoroshev theorems \cite{dumas2014kam}.
They also play a central role in the quantization of classical Hamiltonian systems (in fact, the action and the constant $\hbar$ have the same units) \cite{stone2005einstein, keller1958corrected, keller1960asymptotic, keller1985semiclassical}. 

This paper is outlined as follows. In Section \ref{sec: Ham} we describe the Hamiltonian system described by the Morse potential \eqref{eq:mpot}. In Section \ref{sec:linstab} we describe the equilibria and determine their stability properties in the linear approximation. In Section \ref{eq:po} we discuss the geometry of the region of bounded motion, i.e. the region of periodic orbits bounded by a homoclinic orbit (''separatrix''). In Section \ref{sec:period} we compute the period of the periodic orbits and show how it depends on the energy and other system parameters. In Section \ref{sec:traj} we compute explicit expressions for the trajectories in the region of bounded motion. In Section \ref{sec:homo} we compute an explicit expression for the homoclinic  orbit. In Section \ref{sec:AA} we compute the action-angle variables. In Section \ref{sec:Dleqh} we compute explicit expressions for the trajectories in the region of unbounded motion. In section \ref{sec:app} we apply our results to exploring the nature of chaotic dynamics in the time periodically forced Morse oscillator  and in Section \ref{sec:concl} we discuss our conclusions.

\section{The Hamiltonian}
\label{sec: Ham}

The dynamical system defined by the Morse potential is Hamiltonian, with Hamiltonian function given by:

\begin{equation}
H (q, p) = \frac{p^2}{2m} + D \left( 1-e^{-\alpha q} \right)^2, \qquad (q, p) \in \mathbb{R}^2,
\label{eq:ham}
\end{equation}

\noindent
and Hamilton's equations defined by:

\begin{eqnarray}
\dot{q}  =  \frac{\partial H}{\partial p} & = &  \frac{p}{m}, \nonumber \\
\dot{p}  = - \frac{\partial H}{\partial q} 
& = & -2D \alpha \left(e^{-\alpha q} - e^{-2\alpha q}  \right). 
\label{eq:hameq}
\end{eqnarray}

\noindent
The level sets of the Hamiltonian function have the form:

\begin{equation}
\left\{(q, p) \in \mathbb{R}^2 \, | \, H(q, p) = h = \mbox{constant} \right\}.
\label{eq:levelset}
\end{equation}

\noindent
They are (in general) one dimensional curves that are invariant under the Hamiltonian dynamics, i.e. the Hamiltonian vector field is tangent to the level sets. Since the trajectories lie on these curves, we can use the form of the level curves to obtain parametrizations of  certain trajectories, as we will demonstrate.

\subsection{Equilibria and their Linearized Stability}
\label{sec:linstab}

It is straightforward to verify the the following two points are equilibria for Hamilton's equations:

\begin{equation}
(q, p) = (\infty, 0), \, (0, 0).
\label{eq:equil}
\end{equation}

Next we check their linearized stability properties. The Jacobian matrix, denoted $J$, of Hamilton's equation is given by:

\begin{equation}
J = \left(
\begin{array}{cc}
0 & \frac{1}{m} \\
- 2D \alpha^2 \left(- e^{-\alpha q} + 2 e^{-2\alpha q}  \right) & 0
\end{array}
\right).
\label{eq:jac}
\end{equation}

\noindent
The eigenvalues of $J$ are given by:

\begin{equation}
\pm \sqrt{\rm{det}\, J}.
\label{eq:eivs}
\end{equation}

\noindent
Hence for the two equilibria we have:

\begin{equation}
(q, p) = (0, 0) \Rightarrow \rm{det} \, J = -\frac{2D \alpha ^2}{m},
\label{eq:eqstab}
\end{equation}

\noindent
with corresponding eigenvalues:

\begin{equation}
\pm i\sqrt{\frac{2D}{m}} \alpha,
\label{eq:imageivs}
\end{equation}

\noindent
and

\begin{equation}
(q, p) = (\infty, 0) \Rightarrow \rm{det} \, J =0,
\label{eq:eqsad}
\end{equation}

\noindent
where both eigenvalues are zero.

The equilibrium $(q, p) = (0, 0)$ is stable (''elliptic'' in the Hamiltonian dynamics terminology) and the
 $(q, p) = (\infty, 0)$ is unstable (a ''parabolic'' saddle point  in the Hamiltonian dynamics terminology).

\section{The Region of Bounded Motion: Periodic Orbits}
\label{eq:po}

Using the Hamiltonian \eqref{eq:ham} it is straightforward to verify that the equilibria have the following energies:

\begin{equation}
(q, p) = (\infty, 0) \Rightarrow H=D >0.
\end{equation}

\noindent
and

\begin{equation}
(q, p) = (0, 0) \Rightarrow H=0.
\end{equation}

\begin{figure}
 \centering
 \includegraphics[width=0.49\textwidth]{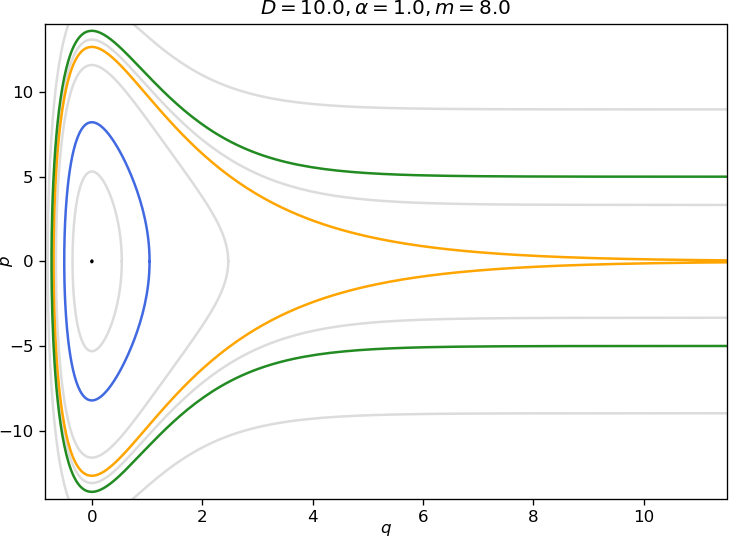}
 \caption{Phase portrait of the Morse oscillator for $D=10$, $\alpha=1$, $m=8$. The equilibrium point at the origin is shown in black, the homoclinic orbit is shown in orange and examples of a periodic orbit and an unbounded trajectory are highlighted with blue and green respectively.}
 \label{fig:phase}
\end{figure}

\noindent
Trajectories with energies larger than $D$  have unbounded motion in $q$. Trajectories having energies $h$ satisfying $0< h < D$ correspond to periodic motions. The level sets of these periodic orbits are given by:

\begin{equation}
h= \frac{p^2}{2m} + D \left( 1-e^{-\alpha q} \right)^2, \quad 0< h < D,
\label{eq:tp1}
\end{equation}

\noindent
and surround the stable equilibrium point $(q, p) = (0, 0)$ as shown in figure \ref{fig:phase}. The periodic orbits intersect the $q$ axis at two distinct points, $q_+ > 0$ and $q_- <0$,  which are referred to as {\em turning points}. These turning points are computed as follows.

Rewriting \eqref{eq:tp1} gives:

\begin{equation}
\frac{p^2}{2m} = h- D \left( 1-e^{-\alpha q} \right)^2.
\label{eq:tp2}
\end{equation}

\noindent
The turning points are obtained from \eqref{eq:tp2} by setting $p=0$:

\begin{equation}
h=D \left( 1-e^{-\alpha q} \right)^2
\label{eq:tp3}
\end{equation}

\noindent
Note that we have:

\begin{equation}
0 \le \sqrt{\frac{h}{D}} \le 1,
\label{eq:tp4}
\end{equation}

\noindent
from which we obtain the following relations:

\begin{eqnarray}
&& 1 \le 1 + \sqrt{\frac{h}{D}} \le 2, \label{eq:tp5} \\
&& 0 \le 1 -  \sqrt{\frac{h}{D}} \le 1, \label{eq:tp6}
\end{eqnarray}

Taking the positive root of \eqref{eq:tp3} gives:

\begin{equation}
1-e^{-\alpha q} =  \sqrt{\frac{h}{D}}.
\label{eq:tp7}
\end{equation}

\noindent
from which we obtain the positive turning point:

\begin{equation}
q_+ = - \frac{1}{\alpha} \log \left(1-  \sqrt{\frac{h}{D}}  \right) > 0.
\label{eq:tppos}
\end{equation}

\noindent
Taking the negative root of \eqref{eq:tp3} gives:

\begin{equation}
1-e^{-\alpha q} = -\sqrt{\frac{h}{D}},
\label{eq:tp8}
\end{equation}

\noindent
from which we obtain the negative turning point:

\begin{equation}
q_- = - \frac{1}{\alpha} \log \left(1 + \sqrt{\frac{h}{D}}  \right) < 0.
\label{eq:tpneg}
\end{equation}

The  level curve with energy equal to the dissociation energy $h=D$ has the form:

\begin{equation}
D= \frac{p^2}{2m} + D \left( 1-e^{-\alpha q} \right)^2,
\label{eq:homo1}
\end{equation}

\noindent
and is a separatrix connecting the (parabolic) saddle point. In the terminology of Hamiltonian dynamics it is a homoclinic orbit. It separates bounded from unbounded motion as illustrated in figure \ref{fig:phase}.

\section{Calculation of the Period of a Periodic Orbit}
\label{sec:period}

In this section we calculate the period of the periodic orbits. 
From \eqref{eq:hameq} we have $\dot{q} = \frac{p}{m}$. Using this expression, and the expression for the level set of the Hamiltonian defining a periodic orbit given in \eqref{eq:tp1}, we have

\begin{equation}
\frac{dq}{dt} = \pm \sqrt{\frac{2}{m}} \sqrt{h-D \left( 1-e^{-\alpha q} \right)^2},
\label{eq:po1}
\end{equation}

\noindent
or

\begin{equation}
\frac{dq}{ \sqrt{h-D \left( 1-e^{-\alpha q} \right)^2}} = \pm \sqrt{\frac{2}{m}} dt.
\label{eq:po2}
\end{equation}

\noindent
We denote the period of a periodic orbit corresponding to the level set with energy value $h$ by $T(h)$. We can obtain the period by integrating $dt$ around this level set.  Using \eqref{eq:po2}, this becomes:

\begin{eqnarray}
T(h) & = & \sqrt{\frac{m}{2}} \int_{q_+}^{q_-} \frac{dq}{ \sqrt{h-D \left( 1-e^{-\alpha q} \right)^2}} -
\sqrt{\frac{m}{2}} \int_{q_-}^{q_+} \frac{dq}{ \sqrt{h-D \left( 1-e^{-\alpha q} \right)^2}}. \nonumber \\
& = & \sqrt{2m}  \int_{q_+}^{q_-} \frac{dq}{ \sqrt{h-D \left( 1-e^{-\alpha q} \right)^2}}.
\label{eq:po3}
\end{eqnarray}

\noindent
Computation of this integral is facilitated by the substitution:

\begin{equation}
u = e^{-\alpha q}.
\label{eq:subs}
\end{equation}

\noindent
After computing the integral using integral 2.266 in \cite{gradshteyn1980table} we obtain:

\begin{equation}
T(h) = \frac{\pi \sqrt{2m}}{\alpha \sqrt{D-h}}
\label{eq:period}
\end{equation}

There are two limits in which this expression can be checked with respect to previously obtained results.
First, we note that for $h=D$ (i.e. the energy of the homoclinic orbit, or ''separatrix'')  $T(D) = \infty$, which is what we expect for the ''period''  of a separatrix.  

Second, we consider $h=0$, which is the energy of the elliptic equilibrium point. In this case we have $T(0) = \frac{\pi \sqrt{2m}}{\alpha \sqrt{D}}$, which is $2 \pi$ divided by the imaginary part of the magnitude of eigenvalue of the Jacobian evaluated at the stable equilibrium point, as we expect.

\section{Expressions for $q(t)$ and $p(t)$, $0 < h < D$}
\label{sec:traj}

In this section we derive an expression for $q(t)$. Differentiating the expression for $q(t)$ will give the expression for $p(t)$ through the relation $\dot{q} = \frac{p}{m}$. Using \eqref{eq:po2}, we have

\begin{equation}
\int_{q_+}^q \frac{dq'}{ \sqrt{h-D \left( 1-e^{-\alpha q'} \right)^2}} =  \sqrt{\frac{2}{m}} t.
\label{eq:traj1}
\end{equation}

\noindent
Choosing  the lower limit if the integral to be $q_+$ is arbitrary, but it is equivalent to the  choice of an initial condition. After computing this integral, we obtain:

\begin{equation}
q(t) = \frac{1}{\alpha} \log{\frac{\sqrt{Dh} \cos \left(\sqrt{\frac{2(D-h)}{m}} \, \alpha t \right) + D}{D-h}}.
\label{eq:q(t)}
\end{equation}

\noindent
It is straightforward to check that the period of \eqref{eq:q(t)} is \eqref{eq:period}.

\subsection{The Homoclinic Orbit}
\label{sec:homo}

As noted earlier, the homoclinic orbit, corresponding to $h=D$, is given by the level set
\eqref{eq:homo1}. Hence, the integral expression for  the homoclinic orbits is obtained from \eqref{eq:po2} by setting $h=D$.  Computing the integral gives:

\begin{equation}
q_0(t) = \frac{1}{\alpha} \log{\frac{1 + \frac{2D}{m} \alpha^2 t^2}{2}}.
\label{eq:qhom}
\end{equation}

\noindent
It is a simple matter to check that $\lim_{t \rightarrow \pm \infty} q(t) = \infty$. Subsequently we obtain $p_0(t)$ from $\dot{q}=\frac{p}{m}$ as

\begin{equation}
p_0(t) = \frac{4mD\alpha t}{2D\alpha^2 t^2+m}.
\label{eq:phom}
\end{equation}

\section{Expressions for Action and the Angle, $0 < h < D$}
\label{sec:AA}

In this section we compute the action-angle representation of the orbits in the bounded region following \cite{melnikov1963vk, wiggins1990introduction, mezic1994integrability}.

We consider a level set defined by the Hamiltonian \eqref{eq:levelset}, for $0 < h < D$. i.e. we consider a periodic orbit with period $T(h)$. Choosing an arbitrary reference point on the periodic orbit, the angular displacement of a trajectory starting from this reference position after time $t$ is given by:

\begin{equation}
\theta = \frac{2 \pi}{T(h)} \int{0}^{t}dt' = \frac{2 \pi}{T(h)}  \sqrt{\frac{m}{2}} \int_{q_+}^q \frac{dq'}{ \sqrt{h-D \left( 1-e^{-\alpha q'} \right)^2}}.
\label{eq:angle1}
\end{equation}

Using the substitution \eqref{eq:subs} and integral 2.266 in \cite{gradshteyn1980table}, we have:

\begin{equation}
\theta = \frac{\pi}{T(h)\alpha} \sqrt{\frac{2m}{D-h}} \left(\frac{3\pi}{2} - \sin^{-1}\frac{(h-D)e^{\alpha q}+D}{\sqrt{Dh}} \right).
\label{eq:angle2}
\end{equation}

\noindent
The action associated with this periodic orbit is the area that it encloses (in phase space) divided by $2 \pi$:

\begin{equation}
I = \frac{1}{2 \pi} \oint_{H=h}  p dq.
\label{eq:action1}
\end{equation}

\noindent
Recalling \eqref{eq:po2}

\begin{equation}
dq = \pm \sqrt{\frac{2}{m}}\sqrt{h-D \left( 1-e^{-\alpha q'} \right)^2} \, dt,
\label{eq:action2}
\end{equation}

\noindent
we obtain:

\begin{equation}
pdq = m \dot{q} dq = \pm \sqrt{2m}\sqrt{h-D \left( 1-e^{-\alpha q'} \right)^2} \, \dot{q}dt,
\label{eq:action3}
\end{equation}

\noindent
and therefore

\begin{equation}
I = \frac{\sqrt{2m}}{\pi} \int_{q_+}^{q_-}\sqrt{h-D \left( 1-e^{-\alpha q'} \right)^2} \, dq'.
\label{eq:action4}
\end{equation}

Using the substitution \eqref{eq:subs} and integral 2.267 in \cite{gradshteyn1980table}, we have:

\begin{equation}
I = \frac{\sqrt{2m}}{\alpha} \left(\sqrt{D}-\sqrt{D-h}\right).
\label{eq:action5}
\end{equation}

\section{Expressions for $q(t)$ and $p(t)$, $0 < D < h$}
\label{sec:Dleqh}

Increasing the total energy $h$ to the value of $D$ and beyond results in unbounded motion. Trajectories retain the turning point $q_-$, while $q_+$ becomes infinite. The expression for $q_-$ is identical to low energies and is obtained from \eqref{eq:ham} by setting $p=0$ and solving for $q$. Recall from \eqref{eq:tpneg} that $q_- = - \frac{1}{\alpha} \log \left(1 + \sqrt{\frac{h}{D}}  \right)$.

For unbounded trajectories it is not possible to define a (finite) period, but we can obtain an expression for $t$ as a function of position. It can be derived by integrating \eqref{eq:po2} from $q_-$ to an arbitrary position $q$ by using the substitution \eqref{eq:subs} as follows:

\begin{eqnarray}
t(q) & = & \sqrt{\frac{m}{2}} \int_{q}^{q_-} \frac{dq'}{ \sqrt{h-D \left( 1-e^{-\alpha q'} \right)^2}}=
       -\frac{1}{\alpha}\sqrt{\frac{m}{2}} \int_{e^{-\alpha q}}^{1+\sqrt{\frac{h}{D}}} \frac{du}{u \sqrt{h-D \left( 1-u \right)^2}},\nonumber \\
       & = & \frac{1}{\alpha}\sqrt{\frac{m}{2(h-D)}} \log{\left(\frac{h-D+De^{-\alpha q}+\sqrt{(h-D)(h-D(1-e^{-\alpha q})^2)}}{\sqrt{hD}e^{-\alpha q}}\right)}.
\label{eq:tinf}
\end{eqnarray}

We obtain an explicit solution $q=q(t)$ by inverting \eqref{eq:tinf}.

\begin{equation}
q(t) = \frac{1}{\alpha} \log{ \frac{\sqrt{hD}e^{2\beta t} -2De^{\beta t} + \sqrt{hD}}{2(h-D)e^{\beta t}} },
\label{eq:qinf}
\end{equation}
where $\beta=\alpha\sqrt{\frac{2(h-D)}{m}}$. Differentiating \eqref{eq:qinf} and using the relation $\dot{q} = \frac{p}{m}$ yields the expression of $p(t)$.

\section{Application: Chaos in the Periodically Forced Morse Oscillator}
\label{sec:app}

Once explicit solutions are known, it is insightful to consider how they change when the vector field is changed, slightly. This is considered in the field of perturbation theory. The subharmonic and homoclinic Melnikov 
\cite{wiggins1990introduction} are examples of global perturbation methods that consider the effect of perturbations on the entire unperturbed integrable system, such as the Morse oscillator that we have considered. The type of perturbation that we will consider is a time periodic excitation of the Morse oscillator. The subharmonic and homoclinic Melnikov methods for this system have been previously considered in \cite{guo2003dynamical}, however, with respect to chaotic dynamics, important technical details for this periodically perturbed Morse oscillator were not considered. We will describe these in detail here, as well as consider the nature of the effect of the parameters of the Morse potential, and the periodic excitation, on chaotic dynamics.

The homoclinic Melnikov function is a well-known and popular method for proving the existence of chaos, in the sense of Smale horseshoes, in time periodically perturbed one degree-of-freedom Hamiltonian systems. The method allows one to prove the existence of transverse homoclinic periodic orbits to a hyperbolic periodic orbit. Such orbits admit a Smale horseshoe construction. There are two issues with the application of the standard homoclinic Melnikov function to the time periodically perturbed Morse oscillator, and  both issues arise from the fact that the saddle point at infinity is parabolic, not hyperbolic. One issue is the fact that the standard derivation of the homoclinic Melnikov function uses hyperbolicity of the fixed point to show that certain terms in the Melnikov function vanish. The other issue involves the construction of the Smale horseshoe map for orbits homoclinic to a parabolic point. Both issue are considered for the perturbed Morse oscillator in \cite{beigie1992dynamics} where it is shown that the growth and decay properties for the parabolic point in the Morse oscillator are sufficient for the standard Melnikov function to be valid as well as the Smale horseshoe map construction to be valid. Hence, the ''Melnikov approach'' is sufficient for determining the existence of chaos for the time-periodically perturbed Morse oscillator. 

We consider the following time-periodic perturbation of the Morse oscillator:

\begin{eqnarray}
\dot{q}  =  \frac{\partial H}{\partial p} & = &  \frac{p}{m}, \nonumber \\
\dot{p}  = - \frac{\partial H}{\partial q} & = & -2D \alpha \left(e^{-\alpha q} - e^{-2\alpha q}  \right) +\varepsilon \cos{\omega t},
\label{eq:perteq}
\end{eqnarray}

\noindent
where $\varepsilon>0$ is the magnitude of the perturbation and $\omega>0$ the frequency of the perturbation. The unperturbed system \eqref{eq:hameq} is obtained by setting $\varepsilon=0$. We formulate the homoclinic Melnikov function using expressions \eqref{eq:qhom}, \eqref{eq:phom} for $q_0(t)$, $p_0(t)$, the homoclinic orbit of the equilibrium point $(q,p)=(\infty,0)$, as follows:

\begin{equation}
    M(t_0,\phi_0)=\int\limits_{-\infty}^{\infty} DH(q_0(t),p_0(t))\cdot(0, \varepsilon \cos(\omega t+\omega t_0+\phi_0)) dt,
\label{eq:Melnikovdef}
\end{equation}

\noindent
where $DH$ is the gradient of the unperturbed Hamiltonian \eqref{eq:ham}, $t_0$ defines the point $(q_0(t_0),p_0(t_0))$ at which the Melnikov funtion $M$ is evaluated and $\phi_0$ is the time at which $M$ is evaluated. We solve \eqref{eq:Melnikovdef} using integral 3.723 in \cite{gradshteyn1980table}.

\begin{eqnarray}
    M(t_0,\phi_0)&=&\varepsilon \int\limits_{-\infty}^{\infty} \frac{\partial H}{\partial p}(q_0(t),p_0(t)) \cos(\omega t+\omega t_0+\phi_0) dt,\nonumber \\
                &=& \varepsilon \int\limits_{-\infty}^{\infty} \frac{p_0(t)}{m} \left( \cos(\omega t)\cos(\omega t_0+\phi_0)-\sin(\omega t)\sin(\omega t_0+\phi_0) \right) dt, \nonumber\\
                &=& -\varepsilon \frac{2m}{\alpha} \sin(\omega t_0+\phi_0) \int\limits_{-\infty}^{\infty} \frac{t}{t^2+\frac{m}{2D\alpha^2}} \sin(\omega t) dt, \nonumber\\
                &=& -\varepsilon \frac{2m}{\alpha} \sin(\omega t_0+\phi_0) e^{-\omega\sqrt{\frac{m}{2D\alpha^2}}}.
\label{eq:Melnikovsol}
\end{eqnarray}

\noindent
Clearly, this function has ''simple zeros'' indicating the existence of homoclinic orbits, and the associated Smale horseshoe type chaos.

The magnitude of the Melnikov function is a measure of the intensity of the chaos. 
Since the exponent in the exponential term is negative, $$e^{-\omega\sqrt{\frac{m}{2D\alpha^2}}}\leq 1.$$ The exponent is decreased by decreasing $\omega$ and $m$, and increasing $\alpha$ and $D$. Note that $M$ is monotonically increasing in dissociation energy $D$, yet the homoclinic orbit ceases to exist for $D=\infty$.

\begin{figure}
 \centering
 \includegraphics[width=0.49\textwidth]{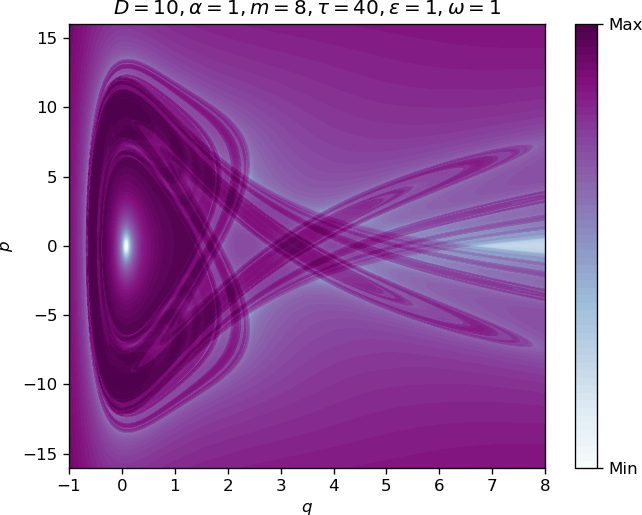}
 \includegraphics[width=0.49\textwidth]{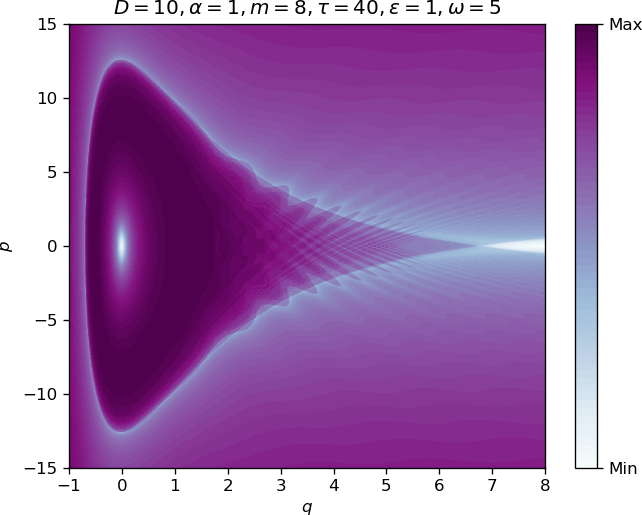}
 \caption{Different intensities of chaos: Poincar{\'e} sections of the system for $t=0$ displaying Lagrangian descriptor values for parameter values $D=10$, $\alpha=1$, $m=8$, $\varepsilon=1$, $\omega=1, 5$. The descriptor is arclength with integration time $\tau=40$. Large chaotic regions (left) for low frequency of perturbation $\omega=1$ and small chaotic regions (right) for high frequency of perturbation $\omega=5$. For sake of clarity, the Lagrangian descriptor values are scaled with arctan. An animation showing the effect of varying $\omega$ can be found at \url{https://youtu.be/6QW7Zjy1poA}. 
 }
 \label{fig:ld}
\end{figure}

The dominant term in the Melnikov functions is $$\frac{2m}{\alpha}.$$ To reach the highest intensity of chaos in the system, $m$ has to be large and $\alpha$ small, whereby the mitigating effect of the exponential term can be countered by choosing $D>\frac{m}{\alpha^2}$. Then for any $\omega$ we can find $(t_0,\phi_0)$, such that $\sin(\omega t_0+\phi_0)=1$ and
\begin{equation}
    |M(t_0,\phi_0)|=\varepsilon \frac{2m}{\alpha} e^{-\omega\sqrt{\frac{m}{2D\alpha^2}}} > \varepsilon \frac{2m}{\alpha}e^{-\frac{\omega}{\sqrt{2}}}.
\label{eq:Melnikovmag}
\end{equation}
In figure \ref{fig:ld} we show a comparison of different intensities of chaos using Lagrangian descriptors, a method introduced by \cite{madrid2009ld} and proven to display invariant manifolds in two dimensional area-preserving maps by \cite{lopesino2015cnsns}.

\section{Conclusion}
\label{sec:concl}

In this paper we have considered the one degree-of-freedom Hamiltonian system defined by the Morse potential energy function (the ''Morse oscillator''). We used the geometry of the level sets to construct explicit expressions for the trajectories as a function of time, their period (for the bounded trajectories, and action-angle variables. We used these trajectories to prove sufficient conditions for chaotic dynamics, in the sense of Smale horseshoes, for the time-periodcally perturbed Morse oscillator using a Melnikov type approach. Using this approach we were able to determine the influence of the parameters of the Morse potential on the chaotic dynamics.

\section*{Acknowledgments}  \noindent We acknowledge the support of  EPSRC Grant no. EP/P021123/1
and   Office of Naval Research (Grant No.~N00014-01-1-0769).

\end{document}